\documentclass[a4paper,12pt]{article}

\usepackage{amssymb,latexsym}
\usepackage{euscript}
\usepackage{eufrak}
\usepackage{amsmath}

\newtheorem{defi}{\bf Definition}[section]
\newtheorem{rem}{\bf Remark}[section]
\newtheorem{lem}{\bf Lemma}[section]
\newtheorem{prop}{\bf Proposition}[section]
\newtheorem{theo}{\bf Theorem}[section]
\newtheorem{cor}{\bf Corollary}[section]
\newtheorem{ex}{\bf Example}[section]

\newenvironment{proof}{\noindent\bf Proof. \rm}{\hfill $\mbox{\boldmath{$ \blacksquare$}}$}

\def\RR{\mbox{{\sl I}}\!\mbox{{\sl R}}} 
\def\fun{\longrightarrow}

\newcommand{\negr}[1]{\boldsymbol{#1}}

\title{\large \bf Monadic distributive lattices and monadic augmented Kripke frames}

\author{A. V. Figallo, I. Pascual and A. Ziliani\\[3mm]
\small Departamento de Matem\'atica, Universidad Nacional del Sur,\\
\small 8000 Bah\'{\i}a Blanca, Argentina\\ 
\small Instituto de Ciencias B\'asicas, Universidad Nacional de San Juan,\\
\small 5400 San Juan, Argentina.
}

\date{}

\begin{document}

\maketitle

\thispagestyle{empty}

\begin{abstract}
In this article, we continue the study of monadic distributive lattices (or m--lattices) which are a natural generalization of monadic Heyting algebras, introduced by Monteiro and  Varsavsky and developed exhaustively by Bezhanishvili. First, we extended the duality obtained by Cignoli for $Q-$distributive lattices to m--lattices. This new duality allows us to describe in a simple  way the subdirectly irreducible algebras in this variety and in particular, to characterize the finite ones.  Next, we introduce the category $\mbox{\boldmath $m{\cal  K F}$}$ whose objects are monadic augmented Kripke frames and whose morphisms are increasing continuous functions verifying certain additional conditions and we prove that it is equivalent to the one obtained above. Finally, we show that the category of perfect augmented Kripke frames given by Bezhanishvili for monadic Heyting algebras is a proper subcategory of $\mbox{\boldmath $m{\cal  K F}$}$.
\end{abstract}

\section{\large Introduction and prelimiaries}\label{S1}

In \cite{FZ}, we introduced monadic distributive lattices as a natural generalization of monadic Heyting algebras defined by Monteiro and Varsavsky in \cite{AM.OV} and deeply studied by Bezhanishvili in \cite{GBz0, GBz, GBz1}.

In \cite{FPZ1}, following the research started in \cite{FZ}, we defined the category $\mbox{\boldmath $m{\cal P}$}$ whose objects are mP--spaces and whose morphisms are mP--functions (\cite[p 7]{FPZ1}) and we showed that it is dually equivalent to the category $\mbox{\boldmath ${\cal M}$}$ of monadic distributive lattices and their corresponding homomorphisms. This duality allowed us to characterized simple and subdirectly irreducible but not simple algebras in this variety. We also proved that the category of topological perfect Ono frames, considered in \cite{GBz} in order to represent monadic Heyting algebras, is a proper subcategory of the one we obtained. Besides, we indicated an example which shows that defining monadic  Heyting  algebras as special distributive lattices they constitute a full subcategory of $\mbox{\boldmath ${\cal M}$}$. 

On the other hand, in \cite{GBz}, Bezhanishvili  established another representation of monadic Heyting  algebras by introducing the category of perfect augmented Kripke frames. The main aim of this paper is to show a new topological representation of monadic distributive lattices in such a way that the latter follows as a particular case of this one. 

\

The paper is organized as follows. In Section 1, we briefly summarize the main definitions and results needed throughout this article. In Section 2, we describe a new topological duality for monadic distributive lattices which extends the one obtained by R. Cignoli for $Q-$distributive lattices. These results enable us to determine a new characterization of the  subdirectly irreducible algebras in this variety.
In Section 3, which is the core of this paper we introduce the notion of monadic augmented   Kripke frame which allows us to consider an equivalent category to the one obtained in  Section 2. Besides, we show that the category of perfect augmented Kripke frames introduced in \cite{GBz} is a proper subcategory of the one we obtained. 

\

In this paper we take for granted the concepts and results on bounded distributive lattices, category theory and universal algebra. To obtain more information on this topics, we direct the reader to the bibliography indicated in \cite{RB.PD}, \cite{SB.HS} and \cite{MacL}. 
However, in order to simplify reading, we summarize the fundamental concepts we use.

Let $X, Y$ be sets. Given a relation $R\subseteq X\times Y$, for each $Z\subseteq X$, $R(Z)$ will denote the image of $Z$ by $R$. If $Z=\{ x\}$, we shall write $R(x)$ instead of $R(\{x \})$. Moreover, for each $V\subseteq Y$, $R^{-1}(V)$ will denote the inverse image of $V$ by $R$, i.e. $R^{-1}(V)= \{ x\in X: R(x)\cap V\not= \emptyset \}$. If $V=\{y\}$, we shall write $R^{-1}(y)$ instead of $R^{-1}(\{y\})$. Besides, if $R,\,T\subseteq X\times X$ then the relation
$R \circ T$ is defined by putting  $(x,y)\in R \circ T$ if and only if there is $z\in X$ such that $(x,z)\in T$ and $(z,y)\in R$.
 
If $X$ is a poset (i.e. partially ordered set) and  $Y \subseteq X $,  then we shall denote by $(Y]$\, ($[Y)$) the set of all $x \in X$ such that $x \leq y$\, ($y \leq x$)
 for some $y \in Y$, and we shall say that $Y$ is increasing (decreasing) if $Y = [Y)$\, ($Y= (Y]$). Furthermore, $max\,Y$ ($min\,Y$) will denote the set of maximal (minimal) elements of $Y$.  

Recall that R. Cignoli (\cite{RC}) introduced the category $\negr{q \cal L}$ of $Q-$distributive lattices and  $Q-$homomorphisms, where a $Q-$distributive lattice is an algebra $\langle L,\vee ,\wedge ,\nabla,0,1 \rangle$ of type $(2,2,1,0,0)$ such that $\langle L, \vee ,\wedge, 0, 1 \rangle$ is a bounded distributive lattice and $\nabla$ is a quantifier on $L$, that is a unary operator on $L$ which satisfies the following equalities: 

\vspace{3mm}

\begin{tabular}{ll}
\hspace{-8mm} (M1) $\nabla 0=0$, & \hspace{2.5cm}(M2) $x\wedge \nabla x=x$,\\[2mm]
\hspace{-8mm} (M3) $\nabla (x\wedge \nabla y)=\nabla x \wedge \nabla y$, & \hspace{2.5cm}(M4)  $\nabla (x\vee y)=\nabla x \vee \nabla y$.\\[2mm] 
\end{tabular}

We shall denote the objects in $\negr{q\cal L}$ by $L$ or $(L,\nabla)$.

\vspace{2mm}

In addition, this author extended Priestley duality (\cite{HP1, HP2, HP3}) to the category $\negr{q\cal L}$. To this aim, he considered  the category  $\negr{q\cal P}$ whose objects are q--spaces and whose morphisms are q--functions. Specifically, a q--space is a pair $(X,E)$ such that $X$ is a Priestley space and $E$ is an equivalence relation on $X$ which satisfies the following conditions:

\begin{itemize}
\item[(E1)] $EU \in D(X)$ for each $U\in D(X)$ where $EU=\{y \in X: (x,y)\in E$ for some $x\in U \}$ and $D(X)$ is the set of all clopen (i.e. simultaneously closed and open) increasing subsets of $X$,

\item[(E2)] the equivalence classes for $E$ are closed in $X$,
\end{itemize}
and a q--function from a q--space $(X,E)$ into another one $(X^{\prime},E^{\prime})$ is an order--preserving continuous function  $f:X \longrightarrow X^{\prime}$ such that $E f^{-1}(U)= f^{-1}(E^{\prime} U)$ for all $U \in D(X^{\prime})$. Besides, he proved:

\begin{itemize}

\item[(E3)] if $L$ is an object in $\negr{q\cal L}$ and $X(L)$ is the Priestley space associated with $L$ (see \cite{HP1, HP2, HP3}), then $(X(L),E_{\nabla})$ is a q--space where $E_{\nabla}$ is the relation defined by $E_{\nabla}=\{(P,R)\in X(L)\times X(L): P\cap \nabla(L)= R\cap \nabla(L)\}$, 

\item[(E4)] if $(X,E)$ is a q--space, then $(D(X),{\nabla}_{E}) $ is a Q--distributive lattice where ${\nabla}_{E}U = EU$ for every $U\in D(X)$.
\end{itemize}

Furthermore, $(L,\nabla)\cong (D(X(L)),{\nabla}_{E_{\nabla}})$ and $(X,E)\cong (X(D(X)),E_{\nabla_E})$ via natural isomorphisms denoted by  $\sigma_{L}$ and $\varepsilon_{X}$ respectively and defined as in \cite{HP1,HP2}. Then, it is concluded that $\negr{q\cal L}$ 
is dually equivalent to $\negr{q{\cal P}}$.

\vspace{3mm}

G. Bezhanishvili (\cite{GBz}) developed the duality theory for monadic Heyting algebras. In order to determine one of these dualities, this author introduced  the category $\mbox{\boldmath$p{\cal AKF}$}$ of perfect augmented Kripke frames (or paK--frames) and their corresponding morphisms (or paK--functions), which we shall des\-cri\-be below. 

A quadruple $(X,\Omega, R,E)$ is a perfect augmented Kripke frame if the following conditions are satisfied:

\begin{itemize}
\item[\rm(k1)] $(X,R,E)$ is an augmented Kripke frame or equivalently

\hspace{-2mm}{\rm(i)} $(X,R)$ is a non--empty partially ordered set,

\hspace{-2mm}{\rm(ii)} $E$ is an equivalence relation on $X$,

\hspace{-2mm}{\rm(iii)} $R\circ E \subseteq E \circ R$.

\item[{\rm(k2)}] $(X,\Omega, R)$ and $(X,\Omega, E\circ R)$ are perfect Kripke frames which means that if  $T=R$ or $T=E\circ R$ then

\hspace{-2mm}{\rm(i)} $(X,T)$ is a non--empty quasi--ordered set,

\hspace{-2mm}{\rm(ii)} $(X,\Omega)$ is a Stone space $($i.e. $0$--dimensional, compact and Hausdorff$)$,

\hspace{-2mm}{\rm(iii)} for all $x \in X$, $T(x)$ is closed in $X$,

\hspace{-2mm}{\rm(iv)} for all clopen $A$ of $X$,  $T^{-1}(A)$ is a clopen of $X$.

\item[{\rm(k3)}] for all increasing clopen $A$ of $X$, $E(A)$ is a clopen of $X$. 
\end{itemize}

Let $(X_1,\Omega_1, R_1,E_1)$ and $(X_2,\Omega_2, R_2,E_2)$ be paK--frames. A paK--function  from $(X_1,\Omega_1, R_1,E_1)$ into $(X_2,\Omega_2, R_2,E_2)$ is a continuous function $f:X_1 \longrightarrow X_2$ which verifies the following conditions: 

\begin{itemize}
\item[\rm (kf1)] for all $x\in X_1$, $y\in X_2$, $(f(x),y)\in T_2$ if and only if there is $z \in X_1$ such that $(x,z)\in T_1$ and $f(z)=y$, where $T_1=R_1$ or $T_1=E_1\circ R_1$, i.e. 
$f$ is strongly isotone with respect to $R_1$ and $E_1\circ R_1$,
  
\item[\rm (kf2)] for all $x\in X_1$, $y\in X_2$,  $(f(x),y)\in E_2$ if and only if there is $z\in X_1$ such that $(x,z)\in E_1$ and $(y, f(z))\in R_2$, i.e. $f$ is almost strongly isotone with respect to $E_1$.
\end{itemize}

Hence, the category of monadic Heyting algebras is dually  equivalent to $\mbox{\boldmath$p{\cal AKF}$}$ (see \cite{GBz}).

\section{\large Subdirectly irreducible algebras in $\mbox{\boldmath ${\cal M}$}$}\label{S2}

Recall that (\cite{FZ}) a monadic distributive lattice {\rm (}or {\rm m}--lattice{\rm)} is a triple $(L,\triangle,\nabla)$ where $L$ is a bounded distributive lattice and $\triangle$, $\nabla$ are unary operations on $L$ satisfying the above mentioned identities {\rm M1}--{\rm M4} and the following ones:

\vspace{3mm}

\begin{tabular}{ll}
\hspace{-8mm} (M5) $\nabla \nabla x= \nabla x$, & \hspace{3.5cm}(M6) $\triangle 1 = 1$,\\[2mm]
\hspace{-8mm} (M7) $x\wedge\triangle x =\triangle x$, & \hspace{3.5cm}(M8)  $\triangle (x \wedge y)= \triangle x \wedge \triangle y$,\\[2mm] 
\hspace{-8mm} (M9) $\triangle \triangle x = \triangle x$, & \hspace{3.5cm}(M10) $\nabla \triangle x = \triangle x$,\\[2mm] 
\hspace{-8mm} (M11) $\triangle \nabla x = \nabla x$. & \\[2mm]
\end{tabular}

The category of m--lattices and their corresponding homomorphisms will be denoted by $\mbox{\boldmath ${\cal M}$}$.

\vspace{3mm}

In what follows firstly, we will determine a topological duality for these algebras which extends the results obtained in \cite{RC}.

\begin{defi}\label{D21}
A q--space $(X,E)$ is an mq--space if the following conditions are fulfilled: 

\begin{itemize}
\item[{\rm(mq1)}] $(x, y)\in E$ and $y \leqslant z$ imply that there is $w \in X$ such that $x \leqslant w$ and $(w, z)\in E$,
\item[{\rm (mq2)}] for every $V\in D(X)$, $(E(X\setminus V)]$ is a clopen subset of $X$.
\end{itemize}
\end{defi}

\begin{defi}\label{D22}
Let $(X_1,E_1)$ and $(X_2,E_2)$ be mq--spaces. An mq--function $f:X_1 \fun X_2$ is a q--function which verifies
\begin{itemize}
\item[{\rm (mqf1)}] $(E_1(f^{-1}(X_2\setminus V))] = f^{-1}((E_2(X_2 \setminus V)])$ for all $V\in D(X_2)$.
\end{itemize}
\end{defi}

We will denote by $\mbox{\boldmath $m{\cal Q}$}$ the category whose objects are mq--spaces and whose morphisms are mq--functions. 

\vspace{2mm}

\begin{rem}\label{R21}
Condition {\rm (mq1)} in Definition \ref{D21} is equivalent to each of the following ones:

\begin{itemize}
\item[{\rm (i)}] $[E(x)) \subseteq E([x))$ for all $x \in X$,
\item[{\rm (ii)}] $E((x]) \subseteq (E(x)]$ for all $x \in X$.
\end{itemize}
\end{rem}

\vspace{2mm}

The properties of mq--spaces and mq--functions that follow are necessary to prove that $\mbox{\boldmath ${\cal M}$}$ and $\mbox{\boldmath $m{\cal Q}$}$ are dually equivalent.

\begin{lem}\label{L21}
Let $(X,E)$ be a q--space. Then $E$ is a closed relation, i.e. $E(A)$ is a closed set for every closed subset $A$ of $X$.
\end{lem}

\begin{proof} 
Let $x \not\in E(A)$. Then, it follows that $(x,y)\not\in E$ for all $y \in A$. Hence, from \cite[Lemma 2.5]{RC}  we conclude that for any $y \in A$ there is a clopen subset $V_y$ of $X$ such that $x \in V_y$, $y \not\in V_y$ and $V_y = E(V_y)$. Therefore,  $A\subseteq \,\bigcup\limits_{y\in A}(X \setminus V_y)$ and using compactness of $A$ we infer that $A \subseteq \bigcup\limits_{i=1}^{n}\ (X \setminus V_{y_i})$.
Since $E(\bigcup\limits_{i =1}^{n} (X \setminus V_{y_i}))= \bigcup\limits_{i =1}^{n} (X \setminus V_{y_i})= X \setminus \bigcap\limits_{i =1}^{n} V_{y_i}$ we have that $\bigcap\limits_{i =1}^{n} V_{y_i} \cap E(A )= \emptyset$. This last assertion and the fact that  $\bigcap\limits_{i =1}^{n} V_{y_i}$ is open and $x \in \bigcap\limits_{i =1}^{n} V_{y_i}$, we deduce that $x \not\in \overline{E(A)}$. Hence, we conclude the proof.
\end{proof}

\begin{cor}\label{C21}
If $(X,E)$ is a q--space, then $E([x))$ is a closed set for all $x \in X$.
\end{cor}

\begin{proof} 
Since $X$ is a Hausdorff space, $\{x\}$ is a closed set for each $x \in X$. Besides, as $X$ is a  Priestley space we have that $[x)$ is closed. Then, by  Lemma \ref{L21} the proof is complete.
\end{proof}

\begin{lem}\label{L22}
If $(X,E)$ is an mq--space, then $E([x))$ is an increasing subset of $X$ for all $x\in X$.
\end{lem}

\begin{proof}
Let $x, y, z \in X$ be such that $y \in E([x))$ and (1) $y \leqslant z$. Then, there is  $w \in X$ such that (2) $x \leqslant w$ and (3) $(y,w)\in E$. From (1), (3) and (mq1), there is  $t \in X$ which verifies (4) $w \leqslant t$ and  (5) $(z,t )\in E$. Hence, by  (2) and (4) we have that $x \leqslant t$. This assertion and (5) imply that $z \in E([x))$. 
\end{proof}

\begin{lem}\label{L23}
Let $(X,E)$ be an mq--space. Then the following con\-di\-tio\-ns are equivalent: 
\begin{itemize}
\item[{\rm (i)}] $x \in X\setminus (E(X\setminus U)]$ for each $U \in D(X)$,
\item[{\rm (ii)}] $E([x))\subseteq U$.
\end{itemize}
\end{lem}

\begin{proof} 
Since $E$ is an equivalence relation, each of the following conditions is equivalent to the next one: (1) $x \in X\setminus (E(X\setminus U)]$,
(2) $[x) \cap E(X\setminus U)= \emptyset$, (3) $E([x))\cap (X\setminus U)= \emptyset$, 
(4) $E([x))\subseteq U$.
\end{proof}

\begin{prop}\label{P21}
If $(X,E)$ is an mq--space, then $(D(X),\triangle _{E},\nabla _{E})$ is an m--lattice, where $\nabla _E(V)= E(V)$ and $\triangle_{E} V= X\setminus (E(X\setminus V)]$ for each $V\in D(X)$.
\end{prop}

\noindent{\bf Proof.}
By virtue of the results established in \cite{RC} it follows that $(D(X),\nabla _{E})$ is a Q--distributive lattice. Besides, from Definition \ref{D21} we have that (M6) to (M9) hold. Then, it only remains to prove that for each $V\in D(X)$ the following conditions are verified:

(i)\, $\nabla _{E}\triangle _{E} V\subseteq \triangle _{E}V$. Let $x\in \nabla _{E}\triangle _{E}V$. Then, there is  (1) $y\in \triangle _{E}V$  such that  $x\in E(y)$. From (1) and Lemma \ref{L23} we have that (2) $E([y))\subseteq V$. On the other hand, by (mq1)we infer that $E([x))= E([y))$ which allows us to conclude by (2) that $E([x))\subseteq V$. This assertion and Lemma \ref{L23} imply that $x\in \triangle _{E}V$.

\vspace{2mm}

(ii)\,$\nabla _{E}V \subseteq \triangle _{E} \nabla _{E}V$. Let $x\in \nabla _E V$. Taking into account that $\nabla _E V$ is an increasing subset of $X$ we infer that $[x)\subseteq \nabla _E V$ and so $E([x)) \subseteq E(V)=\nabla _E V$. Then, by Lemma \ref{L23} we have that $x\in \triangle _{E} \nabla _{E}V$. \hfill $\blacksquare$                                                                                                                                                                                     

\
 
Remark \ref{R22} and Lemmas \ref{L24} and \ref{L25} will be fundamental in order to prove Proposition \ref{P22}.

\vspace{2mm}

\begin{rem}\label{R22}
Let $(L,\triangle, \nabla)$ be an m--lattice and $P,\, T$ filters of $L$. Then, $\triangle^{-1}(T)\subseteq P$ if and only if\, $T \cap \triangle(L)\subseteq P$.
\end{rem}

\begin{lem}\label{L24}{\rm(}see {\rm\cite{FPZ1})}
Let $L$ be a distributive lattice, $\triangle$ an interior operator on $L$ and $a\in L$. If $F\subseteq L$ is a filter such that $\triangle a \notin F$, then there is $Q \in X(L)$ satisfying the following conditions:

\begin{itemize}
\item[{\rm (i)}] $a\notin Q$,
\item[{\rm (ii)}] $\triangle^{-1}(F) \subseteq Q $.
\end{itemize}
\end{lem}

\begin{lem}\label{L25}{\rm(}see {\rm\cite{FPZ1})}
Let $L$ be a distributive lattice and $\triangle$ an interior o\-pe\-ra\-tor on $L$. 
If $S,\,T\in X(L)$ are such that $T\cap \triangle(L)\subseteq S$ then there is $R\in X(L)$ satisfying the following conditions:

\begin{itemize}
\item[{\rm(i)}] $T\subseteq R$,
\item[{\rm(ii)}] $R\cap \triangle(L)=S\cap \triangle(L)$.
\end{itemize}
\end{lem}

\begin{prop}\label{P22}
If $(L,\triangle, \nabla)$ is an m--lattice, then $(X(L), E_{\nabla})$ is an mq--space,  where ${E_\nabla} = \{(P, T)\in X(L)\times X(L): P \cap \nabla(L) = T \cap \nabla(L)\}$. Besides, $\sigma_{L}$ is an m--isomorphism.
\end{prop}

\begin{proof}
From the results established in \cite{RC} we only have to prove that conditions (mq1) and (mq2) in  Definition \ref{D21} hold.

\vspace{2mm}

(mq1)\, Let $(P,Q)\in E_\nabla$ and suppose that $Q \subseteq R$. Then, it follows that $P \cap \nabla(L)\subseteq R \cap \nabla(L)$ and by \cite[Theorem 2.2]{RC} we have that there is  $S\in X(L)$ such that $P\subseteq S$ and $(S,R)\in E_{\nabla}$.

\vspace{2mm}

(mq2)\, Taking into account that for each $V\in D(X)$ there is $a \in L$ such that $V=\sigma_{L}(a)$, it is enough to prove that $( E_\nabla( X(L)\setminus \sigma_{L}(a))]= X(L) \setminus \sigma_{L}(\triangle a)$. Let $P\in (E_\nabla(X(L)\setminus \sigma_{L}(a))]$, then 
$[P) \cap E_\nabla(X(L)\setminus \sigma_{L}(a))\not=\emptyset$.  Therefore, there is $H\in X(L)$ such that $H \in  E_{\nabla}([P))$ and  $a\notin H$. This assertion implies that $\triangle a \notin P$ and so, $P\in X(L)\setminus\sigma_{L}(\triangle a)$. Suppose now that $P\in X(L)$ and $\triangle a\notin P$. Then, by Lemmas \ref{L24} and \ref{L25} there is $Q\in X(L)$ such that $a\notin Q$ and $Q\in E_\nabla([P))$. Therefore,\,  $Q\in X(L)\setminus \sigma_{L}( a)$ and $[P)\cap E_{\nabla}(X(L)\setminus \sigma_{L}(a)) \not=\emptyset$. Hence, we conclude that \, $P\in (E_\nabla(X(L)\setminus \sigma_{L}(a))]$. 
\end{proof}

\vspace{2mm}

Next, our attention is focused on studying some properties of order-preserving continuous functions between mq--spaces, in order to obtain another description of mq--functions more convenient for our purpose.

\begin{lem}\label{L26}
Let $(X_1,E_1)$ and $(X_2, E_2)$ be q--spaces. If $f:X_1 \fun X_2$ is a q--function then the following condition holds:

\begin{itemize}
\item[{\rm (qf1)}] $(E_1(f^{-1}(X_2\setminus V))] \subseteq  f^{-1}((E_2(X_2 \setminus V)])$ for each $V\in D(X_2)$.
\end{itemize}
\end{lem}

\begin{proof} 
Assume that $z\in (E_1(f^{-1}(X_2\setminus V))]$. This means that there are $v, \, u \in X_1$ such that (1) $z \leqslant v$,  (2) $(v,u)\in E_1$ and (3) $f(u)\in X_{2}\setminus V$. From (1), taking into account that $f$ is an increasing  funciton, it follows that (4) $f(z) \leqslant f(v)$. On the other hand, by (2) and the fact that $f$ is a q--function we have that (5) $(f(v),f(u))\in E_2$. Then, by (3), (4) and (5) we conclude that $z \in  f^{-1}((E_2(X_2 \setminus V)])$.
\end{proof}

\begin{prop}\label{P23}
Let $(X_1, E_1)$ and $(X_2, E_2)$ be mq--spaces. Then for each order--preserving continuous  function $f:X_1 \fun X_2$ the following conditions are equivalent:
\begin{itemize}

\item[\rm (qf1)] $(E_1(f^{-1}(X_2\setminus V)) ] \subseteq  f^{-1}( (E_2(X_2 \setminus V)] )$ for each $V\in D(X_2 )$,
\item[\rm (qf2)] $f(E_1([x))) \subseteq E_2([f(x)))$ for each $x \in X_1$,
\item[{\rm (qf3)}] $[f(E_1([x)))\subseteq  E_2([f(x)))$ for each $x \in X_1$.
\end{itemize}
\end{prop}

\begin{proof} 
(qf1)\, $\Rightarrow$\, (qf2): Let (1) $y \in E_1([x))$ and suppose that (2) $f(y)\notin E_2([f(x)))$. This assertion and Lemma \ref{L22} imply that $z\not\leq f(y)$ for all $z\in E_2([f(x)))$. Then, there is $U_z\in D(X_2)$ such that $z \in U_z$ and $f(y)\in X_2\setminus U_z$ for all $z\in E_2([f(x)))$. Taking into account Corollary \ref{C21}, a simple compactness argument shows that there is $U\in D(X_2)$ such that (3) $E_2([f(x))) \subseteq U$ and (4) $f(y)\in X_2\setminus U$. From (1) and (3) we infer that $E_1([x)) \cap f^{-1}(X_2\setminus U)\not= \emptyset $ and so, $[x) \cap E_1(f^{-1}(X_2\setminus U))\not= \emptyset$. Hence, $x\in (E_1(f^{-1}(X_2\setminus U))]$ and by (qf1) we have that $f(x)\in (E_2(X_2 \setminus U)]$. Therefore, $[f(x))\cap E_2(X_2\setminus U)\not=\emptyset$. this last condition implies that  $E_2([f(x)))\cap (X_2\setminus U)\not=\emptyset$,  which contradicts (3).

\vspace{2mm}

(qf2)\, $\Rightarrow$\, (qf3): It is a direct consequence of the hypothesis and Lemma \ref{L22}.

\vspace{2mm}

(qf3)\, $\Rightarrow$\, (qf1): Let $z\in (E_1(f^{-1}(X_2\setminus V))]$. Then, there are $v, u \in X_1$ such that (1) $z \leqslant v$,  (2) $(v,u)\in E_1$ and (3) $f(u)\in X_{2}\setminus V$. From (1) and (2), $u \in E_1([z))$ and so, $f(u)\in f(E_1([z)))$ which implies by (qf3) that $f(u) \in E_2([f(z)))$. Therefore, there is $t \in X_2$ such that $f(z)\leqslant t$ and   $(f(u),t)\in E_2$. From these last assertions and (3) we conclude that $z \in f^{-1}((E_2(X_2 \setminus V)])$.
\end{proof}

\begin{prop}\label{P24}
Let $(X_1,E_1)$ and $(X_2,E_2)$ be mq--spaces. Then for each order--preserving continuous  function $f:X_1 \fun X_2$ the following conditions are equivalent:

\begin{itemize}
\item[\rm (mqf2)] $f^{-1}((E_2(X_2 \setminus V)])\subseteq (E_1(f^{-1}(X_2\setminus V))]$ for each $V\in D(X_2 )$,
\item[\rm (mqf3)] $E_2([f(x)))\subseteq [f(E_1([x)))$ for each $x \in X_1$.
\end{itemize}
\end{prop}

\begin{proof}
(mqf2)\, $\Rightarrow$\, (mqf3): Let (1) $z \in E_2([f(x))$ and suppose that $f(x)\leqslant z$. Since $f(x)\in f(E_1([x)))$ we conclude that $z \in [f(E_1([x)))$ which completes the proof. Now, suppose that (1) does not hold, then there is $U\in D(X_2)$ such that (2) $z\notin U$ and $f(x)\in U$ from which by (1) we conclude that $x \in f^{-1}(E_2(X_2 \setminus U)]$. Then, by (mqf2) we have that $x\in (E_1(X_1 \setminus  f^{-1}(U))]$. From this assertion it follows that $E_1([x)) \cap (X_1 \setminus  f^{-1}(U))\not= \emptyset$. Furthermore, since $f$ is a continuous and closed function, by Corollary \ref{C21} we deduce that  $K=f(E_1([x))) \cap (X_1 \setminus f^{-1}(U)))$ is a closed subset of $X_2$ and therefore it is compact. Then, there is $y \in E_1([x))\cap (X \setminus f^{-1}(U))$ such that $f(y)\leqslant z$. Indeed, if we suppose that $f(y)\not\leqslant z$ for all $y \in E_1([x))\cap (X \setminus f^{-1}(U))$ then, $k \not\leqslant z$ for all $k \in K$. Hence, a compactness argument shows that there are subsets $V_i \in D(X_2)$,  $1\leqslant i \leqslant n$ such that $z\notin V_i$ and $K \subseteq \bigcup\limits_{i=1}^{n} V_i$. If $V=\bigcup\limits_{i=1}^{n} V_i$ then (3) $K \subseteq V$ and (4) $z \not\in V$. On the other hand, let $W= U \cup V$ and so, by (2) and (4) it follows that $z \not\in W$. This statement and (1) allows us to conclude that  $x \in  f^{-1}( E_2( X_2 \setminus W)])$, from which by (mqf2) we obtain that $x \in (E_1(X_1 \setminus f^{-1}(W)]$. Therefore, $(X_1 \setminus (f^{-1}(U )\cup f^{-1}(V ))\cap E_1([x))\not=\emptyset$ and then $f(E_1([x)) \cap (X_1 \setminus f^{-1}(U)) \cap (X_1 \setminus f^{-1}(V)))\not=\emptyset$. This last assertion implies that $K \cap (X_2 \setminus V) \not=\emptyset$ and then $K \not\subseteq V$, which contradicts (3). Thus, we conclude that there is $y \in E_1([x))$ such that $f(y)\leqslant z$ and so it follows (mqf3).    

\vspace{2mm}

(mqf3)\, $\Rightarrow$\, (mqf2): Let $f(x)\in (E_2(X_2 \setminus V)]$. Then, $E_2([f(x))) \cap (X_2 \setminus V)\not= \emptyset$ and by (mqf3) we have $[f(E_1([x)))) \cap (X_2 \setminus V)\not= \emptyset$. Besides, taking into account that $X_2 \setminus V$ is a decreasing subset of $X_2$ we conclude that $f(E_1([x)))\cap (X_2 \setminus V) \not= \emptyset$. This last assertion implies that $E_1([x))\cap f^{-1}(X_2 \setminus V)\not= \emptyset$ and therefore, $[x)\cap E_1(f^{-1}(X_2 \setminus V))\not= \emptyset$ from which we infer that $x \in (E_1(f^{-1}( X_2 \setminus V))]$.
\end{proof}

\vspace{2mm}

The above results allow us to obtain the description of  mq--functions we were looking for.   

\begin{cor}\label{C22}
Let $(X_1,E_1)$ and $(X_2,E_2)$ be mq--spaces. Then for each order--preserving continuous  function $f:X_1 \fun X_2$\, the following conditions are equivalent:

\begin{itemize}
\item[\rm(i)] $f$ is an mq--function, 
\item[\rm(ii)] $f$ is a q--function which verifies {\rm(mqf3)},
\item[\rm(iii)] $f$ is a q--function which verifies {\rm(mqf4)}: $E_2([f(x)))= [f(E_1([x))))$.
\end{itemize}
\end{cor}

\begin{proof}
(i)\, $\Leftrightarrow$\, (ii): It is a direct consequence of Definition \ref{D22} and Lemma   \ref{L26}.

\vspace{2mm}

(i)\, $\Leftrightarrow$\, (iii): It follows from Definition \ref{D22} and Propositions \ref{P23} and \ref{P24}. 
\end{proof}

\

Next, we are going to characterize the isomorphisms in the category $\mbox{\boldmath $m{\cal Q}$}$ for which Lemma \ref{L27} and Corollary \ref{C23} are fundamental.

\begin{lem}\label{L27}
Let $(X_1,E_1)$ and $(X_2,E_2)$ be mq--spaces. If $f:X_1 \fun X_2$ is an isomorphism in $\negr{q{\cal P}}$, then $ E_2([f(x)))= f(E_1([x)))$ for each $x \in X_1$.
\end{lem}

\begin{proof}
From Lemma \ref{L26} and Proposition \ref{P23} we have that $f(E_1([x))) \subseteq  E_2([f(x)))$ for all $x \in X_1$ and since $ f^{-1}$ is a q--function, we infer also that $f^{-1}(E_2([y))) \subseteq  E_1([f^{-1}(y)))$ for all $y \in X_2$. This statement and the fact that  $f$ is bijective imply that $ E_2([f(x))) \subseteq  f(E_1[x))$ for all $x \in X_1$.
\end{proof}

\begin{cor}\label{C23} 
Let $(X_1, E_1)$ and $(X_2, E_2)$ be mq--spaces. If $f:X_1 \fun X_2$ is an isomorphism in  $\negr{q{\cal P}}$ then, $f$ is an mq--function.  
\end{cor}

\begin{proof}
Note first that (1) $f(E_1([x)))$ is an increasing subset of $X_2$. Indeed, let $y \in f(E_1([x)))$ and (2) $y \leqslant z$. Then, there is (3) $t \in E_1([x))$ such that (4) $f(t)= y$. On the other hand, since $f$ is bijective there is $w \in X_1$ such that (5) $z =f(w)$. Then, taking into account that $f$ is an order isomorphism from (2), (4) and (5) it follows that $t\leqslant w$. This last assertion, (3) and  Lemma \ref{L22} imply that $w \in E_1([x))$. Therefore, by (5) we conclude that  $z \in f(E_1([x)))$. Hence, from (1) we have that  $f(E_1([x)))= [f(E_1([x)))$  for all $x \in X_1$. Then, by Lemma \ref{L27} we obtain that $E_2([f(x)))= [f(E_1([x))))$ for all $x \in X_1$ and hence, by Corollary \ref{C22} we conclude the proof. 
\end{proof}

\begin{prop}\label{P25}
Let $(X_1,E_1)$, $(X_2, E_2)$ be mq--spaces and let $f:X_1 \fun X_2$ be a function. Then the following   conditions are equivalent:

\begin{itemize}
\item[\rm (i)] $f$ is an isomorphism in $\negr{q{\cal P}}$,
\item[\rm (ii)] $f$ is an isomorphism in $\mbox{\boldmath $m{\cal Q}$}$.
\end{itemize}
\end{prop}

\begin{proof} 
Since $f$ and $f^{-1}$ are isomorphisms in $\negr{q{\cal P}}$, by Corollary \ref{C23} they are both  mq--functions and so, we conclude that $f$ is an isomorphism in $\mbox{\boldmath $m{\cal Q}$}$. The proof of the other implication follows immediately from Definition \ref{D22}. 
\end{proof}

\begin{cor}\label{C24}
Let $(X,E)$ be an mq--space. Then $\varepsilon_{X}: X\fun X(D(X))$ is an isomorphism in   $\mbox{\boldmath $m{\cal Q}$}$.
\end{cor}

\begin{proof} 
It is a direct consequence of \cite[Theorem 2.6]{RC} and Proposition \ref{P25}. 
\end{proof}

\

From the above results and using the usual procedures we conclude 

\begin{theo}\label{T21} 
The category $\mbox{\boldmath $m{\cal Q}$}$ is naturally equivalent to the dual of the category    $\mbox{\boldmath ${\cal M}$}$.
\end{theo}

Next, we obtain a characterization of the congruence lattice on  monadic distributive lattices by means of certain closed subsets of its associated mq--space. This fact allows us to describe the congruence lattice on $Q-$distributive lattices completing the results obtained in \cite{RC}. 

\begin{defi}\label{D23} 
Let $(X,E)$ be an  mq--space. A subset $Y$ of $X$ is $id-$saturated provided that 
$min \,E([y)) \cup max \,(E(y))\subseteq  Y$ for all $y \in Y$.
\end{defi}

\begin{theo}\label{T22}
Let $L\in \mbox{\boldmath ${\cal M}$}$ and let $X(L)$ be the mq--space associated with $L$. Then, the lattice ${\cal C}_{idS}(X(L))$ of closed and $id-$saturated subsets of $X(L)$ is isomorphic to the dual lattice $Con_{M}(L)$ of m--congruences on $L$ and the i\-so\-mor\-phism is the function $\Theta_{M}$ defined by the prescription $\Theta_{M}=\{(a,b)\in L\times L: \sigma_L(a) \cap Y=\sigma_L(b) \cap Y\}$. 
\end{theo}

\begin{proof}
Let $Y \in {\cal C}_{idS}(X(L))$. Then by the results established in \cite {HP1}, (see also \cite[Section 6]{HP2} and \cite{HP3}) we have that $\Theta_{M}(Y)$ is a lattice congruence. Hence, to prove that $\Theta_{M}(Y)\in Con_{M}(L)$, we need to show that $\Theta_{M}(Y)$ preserves $\nabla$ and $\triangle$. Let $(a, b) \in \Theta_M(Y)$ and take $P \in \sigma_{L} (\nabla a) \cap Y= \nabla_{{E}_{\nabla} } \sigma_L(a) \cap Y$. Then, we have  that there is $Q \in E_{\nabla}(P)$ such that  $a \in Q$.  Since $(X(L), E)$ is a qP--space, $E_{\nabla}(P)$ is a closed subset of $X(L)$.  Then, there is  $S \in max\,  E_{\nabla}(P)$ such that $Q \subseteq S$. Since $Y$ is $id-$saturated and $P \in Y$ we obtain that $S \in Y$. On the other hand, it follows that $a \in S$. Hence, by the hypothesis we have that $S \in \sigma_L(a) \cap Y=\sigma_L(b) \cap Y$. This assertion implies that $\nabla b \in S$.  Since  $S \cap \nabla(L)= P \cap \nabla(L)$, then  
 $P \in \sigma_{L} (\nabla b) \cap Y$. Therefore, $\sigma_{L}(\nabla a) \cap Y \subseteq \sigma_{L}(\nabla b) \cap Y$. The proof of the other inclusion is similar.   

\vspace{2mm}
    
On the other hand, if we suppose  that $(\triangle a, \triangle b) \notin \Theta_M(Y)$ we can assume without loss of generality that there is $P\in Y$ such that  $\triangle a \in P$ and 
$\triangle b \notin P$. From this last assertion and Lemmas \ref{L24} and \ref{L25} we have that  there is $S \in X(L)$ such that $b \notin S$ and $S \in E_{\nabla}([P))$. Then, by Lemma \ref{L21} we infer that  there is $T \in min \, E_{\nabla}([P))$ such that  $T \subseteq S$. Taking into account that $Y$ is $id-$saturated we have that $T \in Y$. Besides, we conclude that $T \notin \sigma_L(b) \cap Y$. Since  $P \cap \nabla(L) \subseteq T \cap\nabla(L)$ and $\triangle a \in P$ it follows that $T \in \sigma_L(a) \cap Y$. Therefore, $(a, b) \not\in \Theta_M(Y)$ which is a contradiction. 

\vspace{2mm}

Conversely, let $\theta\in Con_{M}(\bf A)$ and let $h\!: L\longrightarrow L/\theta$ be the natural epimorphism. Since $\theta$ is a lattice congruence on $L$, we have that 
$Y=\{h^{-1}(R): R\in X(L/\theta\}$ is a closed subset of $X(L)$ and $\theta=\Theta(Y)$ (see \cite{HP1,HP2,HP3}). Then, it only remains to prove that $Y$ is $id-$saturated. More precisely, 

\vspace{2mm}

(i)\, $min \, \, E_{\nabla}([P)) \subseteq Y$ for each $P \in Y$: Suppose that this statement is not true, then there is  $T\in min  \, \, E_{\nabla}([P))$ and $T\notin Y$. Therefore, there are $a,\, b\in L$ such that $T \in \sigma_L(a)\setminus \sigma_L(b) \subseteq X(L) \setminus Y$. This last assertion implies that $(b,a\vee b)\in \Theta (Y)$. On the other hand, we have that   $P \cap \triangle(L)\subseteq  T \cap \triangle(L)$ and taking into account that $b\notin T$, it follows  that $\triangle b\notin P$. Hence, we conclude that $\triangle (a\vee b)\notin P$. Let us consider the filter $\triangle^{-1}(P)$, the ideal $((L \setminus T) \cup \{a\}]$ and suppose that $((L \setminus T) \cup \{a\}] \cap \triangle^{-1}(P) = \emptyset$. Then, by the Birkhoff--Stone theorem there is $S\in X(L)$ such that \linebreak  $S \subseteq L \setminus ((L \setminus T) \cup \{a\}] \subset T$ and  by Remark \ref{R22}  and  Lemma \ref{L25} it follows that $S\in E_{\nabla}([P))$. These  assertions contradict the fact that  $T \in  min  \, \, E_{\nabla}([P))$. Therefore, $((L \setminus T) \cup \{a\}] \cap \triangle^{-1}(P) \not= \emptyset$. This last statement implies that there is $x \in ((L \setminus T) \cup \{a\}]$ and $\triangle x \in P$. Then, there is $q\in L \setminus T$ such that $\triangle (q \vee a)\in P$ and so,
$\triangle (q\vee a \vee b)\in P$. Since
 $(q\vee b, q\vee a \vee b)\in \Theta(Y)$, we have that $(\triangle (q\vee b), \triangle (q\vee a \vee b))\in \Theta(Y)$. Hence, $q\vee b\in \triangle^{-1}(P)$ and since $\triangle^{-1}(P)\subseteq T$ we obtain that  $q\in T$ or $b\in T$ which is a contradiction. Therefore, $T\in Y$. 

\vspace{2mm}

(ii)\, $max  \, \, E_{\nabla} (P) \subseteq Y$ for all $P \in Y$: Let $P \in Y$ and $T  \in max \,  E_{\nabla} (P)$. Suppose that $T \not\in Y$.  Since $Y$ is a closed subset of $X(L)$, then there are $a, \,b \in L$ such that $T \in \sigma_L(a) \setminus \sigma_L(b)$  and $(\sigma_L (a) \setminus \sigma_L (b)) \cap Y = \emptyset$. This last assertion implies that $(a, a\wedge  b)\in \Theta(Y)$. Since $a \wedge b \notin T$ we infer that $T\subset F$, where $F$ is the filter of $L$ generated by $T\cup \{a \wedge b \}$. On the other hand, from the hypothesis,
$T\cap \nabla L=P\cap \nabla L$ which implies that $P\cap \nabla L\subset F\cap \nabla L$. This means that there is $q \in T$ such that $\nabla(a \wedge b \wedge q)\not\in P$. Since $(a \wedge q, a \wedge b\wedge q)\in \Theta(Y)$ we infer that $(\nabla(a \wedge q),\nabla( a \wedge b\wedge q))\in \Theta(Y)$.\linebreak  Furthermore, $a\wedge q \in T$ and so, we have that 
$\nabla(a\wedge q)\in P$. From this assertion and taking into account that $P\in Y$ we conclude that $\nabla( a \wedge b\wedge q)\in P$ which is a contradiction. 
\end{proof}

\

Then, as a direct consequence of the proof of Theorem \ref{T22}, we obtain Corollary \ref{C25} which determines the congruences on $Q-$dis\-tri\-bu\-ti\-ve lattices.

\begin{defi}\label{D24} 
Let $(X,E)$ be a q--space. A subset $Y$ of $X$ is $i-$saturated provided that 
$max \,E(y)\subseteq  Y$ for all $y \in Y$.
\end{defi}

\begin{cor}\label{C25}
Let $L$ be a $Q-$distributive lattice and let $X(L)$ be the q--space associated with $L$. Then, the lattice  ${\cal C}_{iS}(X(L))$ of closed and $i-$saturated subsets of $X(L)$ is isomorphic to the dual lattice  $Con_{Q}(L)$ of Q--congruences on $L$. 
\end{cor}

Our next task will be to determine the subdirectly irreducible members of $\mbox{\boldmath ${\cal M}$}$ from which Theorem \ref{T22} is fundamental. 

\begin{cor}\label{C26}
Let $(L, \triangle,\nabla)$ be am m--lattice and let $(X(L), E_\nabla)$ be the mq--space associated with $L$. If $Y$ is an $id-$saturated subset of $X(L)$, then $\overline{Y}$ is also $id-$saturated {\rm(}where $\overline{Y}$ denotes the closure of $Y${\rm)}. 
\end{cor}

\begin{proof}
Taking into account that for each $x \in L$, $\sigma_L(x)$ is a clopen subset of $X(L)$, it is easy to check that $(\sigma_L(a)\bigtriangleup \sigma_L(b))\cap \overline Y= \emptyset$ if and only if $(\sigma_L(a)\bigtriangleup \sigma_L(b))\cap \ Y= \emptyset$. Therefore, we infer that $\Theta(Y) =\{(a,b)\! \in \! L\times L: \, \, (\sigma_L(a)\bigtriangleup \sigma_L(b)) \cap Y = \emptyset \} = \{(a,b)\! \in \! L\times L: \, \, \sigma_L(a) \cap \overline Y =  \sigma_L(b) \cap \overline Y \}=\Theta(\overline Y)$. 
Since $Y$ is an $id-$saturated subset of $X(L)$ arguing as in Theorem \ref{T22} we infer that for all  $(a,b) \in \Theta(Y)$ we have that $(\sigma_L(\triangle a)\bigtriangleup \sigma_L(\triangle b)) \cap Y = \emptyset$ and \linebreak $(\sigma_L(\nabla a)\bigtriangleup \sigma_L(\nabla b)) \cap Y = \emptyset $.  Therefore, $(\triangle a, \triangle  b)\in \Theta(\overline Y)$ and $(\nabla a, \nabla  b)\in \Theta(\overline Y)$ for all $(a,b) \in \Theta(\overline Y)$. These assertions imply that $\Theta(\overline Y)\in Con_{M}(L)$. Hence, by Theorem \ref{T22} we conclude that $\overline Y$ is an $id-$saturated subset of $X(L)$.
\end{proof}

\

On each bounded lattice we can define a special quantifier, namely the indiscrete or simple quantifier given by the prescription $\nabla 0=0$ and $\nabla x=1$ for each $x\in L$, $x\not=0$. 
If $L$ is an m--lattice and $\nabla$ is the simple quantifier, taking into account the results established in \cite{FPZ1}, we can assert that $\triangle 1=1$ and $\triangle x=0$ for all $x\in L$, $x\not=1$. In this case, we say that $(\nabla, \triangle)$ is simple and it will play an important role in the characterization of simple m--lattices.

\begin{prop}\label{P26}
Let $(X,E)$ be an mq--space. Then the following conditions are equivalent: 

\begin{itemize}
\item[{\rm(i)}] $\nabla_{E}$ is the simple quantifier,  
\item[{\rm(ii)}] $E= X \times X$.
\end{itemize}
\end{prop}

\begin{proof}
(i)\, $\Rightarrow$\, (ii): Suppose that $E \not= X \times X$. Then, there are $x, y \in X$ such that $(x, y)\notin E$. Since $(X, E)$ is a qP--space from \cite[Lemma 2.5]{RC} there is  $U \in D(X)$, $U = \nabla_{E} U$ such that $x \in U$  and $y \not\in U$  or  $y \in U$ and $x \not\in U$. Hence, $U \not= \emptyset$ and  $\nabla_{E} (U)\not= X$ which contradicts the hypothesis.  

\vspace{2mm}

(ii)\, $\Rightarrow$\, (i): Let $U \in D(X)\setminus \{\emptyset\}$. Since $E(x)=X$ for all $x \in X$, then $\nabla_{E}(U)= X$ which completes the proof.
\end{proof}

\

As a direct consequence we get

\begin{cor}\label{C27}
Let $(X,E)$ be an  mq--space. If $\nabla_E$ is the simple quantifier then $min \, X \cup max\, X$ is $id-$saturated.
\end{cor}

\begin{prop}\label{P27}
Let $(X,E)$ be an  mq--space such that $\nabla_{E}$ is the simple quantifier. Then for all non--empty subset $Y$ of $X$ the following conditions are equivalent:

\begin{itemize}
\item[{\rm(i)}] $Y$ is $id-$saturated,
\item[{\rm(ii)}] $min \, X \cup max\, X \subseteq Y$.
\end{itemize}
\end{prop}

\begin{proof}
It follows immediately from Proposici\'on \ref{P26}.
\end{proof}

\

Corollary \ref{C28} and Lemma \ref{L28} combine are necessary in order to prove Theorem \ref{T23}.

\vspace{2mm}
 
\begin{cor}\label{C28}
Let $(X,E)$ be an mq--space  such that $\nabla_{E}$ is the simple quantifier. Then $\overline{min \, X \cup max\, X}$ is the lowest non--empty closed and $id-$sa\-tu\-ra\-ted subset of $X$.
\end{cor}

\begin{proof}
Since $X$ is a Priestley space it follows that $min\, X \cup max\,  X \not= \emptyset$. Then, from Corollary \ref{C27} and Corollary \ref{C26} we have that $\overline{min \, X \cup max\, X}$ is a non--empty closed and $id-$saturated subset of $X$. On the other hand, if $Y$ is a non--empty closed  and $id-$saturated subset of $X$, by Proposition \ref{P27} we conclude the proof.
\end{proof}

\begin{lem}\label{L28}
Let $(X,E)$ be an mq--space. Then the following conditions hold:

\begin{itemize}
\item[{\rm(i)}] for each $x\in X$, $E([x))$ is a closed and $id-$saturated subset of $X$,
\item[{\rm(ii)}] $\nabla_{E}(U)$ is an $id-$saturated subset of $X$ for all $U \in D(X)$.
\end{itemize}
\end{lem}

\begin{proof}
(i)\, It follows from Corollary \ref{C21}, (mq1) and the fact that $E \circ \leqslant$ is a quasi--order.

\vspace{2mm}

(ii)\, By (M11) we have that $\nabla_{E} U= X \setminus (E( X \setminus \nabla_{E} U)]$ and so, by Lemma \ref{L23} we conclude the proof.
\end{proof}

\begin{theo}\label{T23}
Let $(X,E)$ be an mq--space. Then the following conditions are equivalent: 

\begin{itemize}
\item[{\rm(i)}] $(D(X),\triangle_{E}, \nabla_{E})$ is a simple monadic distributive lattice,
\item[{\rm(ii)}]  $(\triangle_{E}, \nabla_{E})$ is simple and $X = \overline{ min \,X \cup max \, X}$.
\end{itemize}
\end{theo}

\begin{proof}
(i)\, $\Rightarrow$\, (ii): Suppose that $\nabla_{E}$ is not the simple quantifier. Then, there is $U \in D(X)\setminus\{\emptyset, X\}$ such that $ \nabla_{E}\,U \not= X$ and so, by (ii) in Lemma \ref{L28} we conclude that $ \nabla_{E}\,U $ is a proper non--empty closed and $id-$saturated subset of $X$.  Therefore, by Theorem \ref{T22} we have that $(D(X),\triangle_{E}, \nabla_{E})$ is not a simple m--lattice which contradicts (i). Hence, the fact that $\nabla_{E}$ is simple, Corollary \ref{C28} and (i) imply that $X = \overline{ min \,X \cup max \, X}$.

\vspace{2mm}

(ii)\, $\Rightarrow$\, (i): From the hypothesis and Corollary \ref{C28} we have that \\ ${\cal C}_{idS}(X(D(X)))=\{\emptyset, X(D(X)) \}$. Hence, by Theorem \ref{T22} the proof is complete.
\end{proof}

\begin{prop}\label{P28} 
Let $(L, \triangle, \nabla)$ be an m--lattice and let $(X(L),E_\nabla)$ be the mq--space associated with $L$. Then the following conditions are equivalent:
    
\begin{itemize}
\item[{\rm(i)}] $\overline{max\, X \cup \, min \, X} =X$,
\item[{\rm(ii)}] for each $a,\,b\in L$ if $ a \not\leqslant b$, there is $M \in max \, X(L)$ such that $a \in M$ and $b \not\in M$ or there is $P \in min \, X(L)$ such that $a \in P$ and $b \not\in P$,
\item[{\rm(iii)}] for each $a,\,b\in L$ if $a \not\leqslant b$ there is $c \in L$ such that $a \wedge c \not= 0$ and $b \wedge c = 0$, or there is $d \in L$ such that $b \vee d \not= 1$ and $a \vee d = 1$.
\end{itemize}
\end{prop}

\begin{proof}
(i)\, $\Leftrightarrow$\, (ii):  It follows from the fact that the non--empty basic open sets of  $X(L)$ are $\sigma_L(a)\setminus \sigma_L(b)$  with $a,\, b \in L$,  $a \not\leqslant b$. 

\vspace{2mm}

(ii)\, $\Rightarrow$\, (iii): Suppose that there is $M \in max \, X(L)$ such that $ a \in M$ and  $b \not\in M$. Then, we have that $M \in \sigma_L(a)$ and $M \notin \sigma_L(b)$. Since $M \in max \, X(L)\setminus \sigma_L(b)$ we conclude that $M \not\subseteq P$ for all $P \in \sigma_L(b)$. This last assertion implies that for all $P \in \sigma_L(b)$ there is  $c_P \in L$ such that $M \in \sigma_L(c_P)$ and $P \not\in \sigma_L(c_P)$. A compactness argument shows that there are $c_{P_1}, \ldots, c_{P_n}\in L$ such that $\sigma_L(b) \cap \bigcap\limits_{i =1}^{n} \sigma_L(c_{P_i})= \emptyset$ and $M \in \bigcap\limits_{i=1}^{n} \sigma_L(c_{P_i})$. Take $c = \bigwedge\limits_{i=1}^{n}c_{P_i}$, then $\sigma_L(b) \cap \sigma_L(c)= \emptyset$ and therefore $b \wedge c = 0$. On the other hand, we have that $c \in M$ and since $a \in M$ we infer that $a \wedge c \not= 0$.

Suppose now that there is $Q \in min \, X(L)$ such that $ a \in Q$ and $b \not\in Q$. Then, arguing as in the previous case we conclude that $a \vee d = 1$ and $b \vee d \not= 1$ for some $d\in L$.
 
\vspace{2mm}

(iii)\, $\Rightarrow$\, (ii): Let  $a,\,b\in L$ be so that $a \not\leqslant b$ and suppose that there is $c \in L$ such that (1) $a \wedge c \not= 0$ and (2) $b \wedge c = 0$. Then, by (1) there is $P \in X(L)$ and (3) $a \wedge c \in P$. Since $X(L)$ is a Priestley space, there is (4) $M \in max \, X(L)$ and (5) $P\subseteq M$. Therefore, $a \wedge c \in M$ which implies that $a \in M$.  On the other hand, from (3) and (5) we have that $c \in M$ and so, by (2) and (4) we conclude that $b \notin M$.

Suppose now that there is $d \in L$ and (6) $b \vee d \not= 1$,  (7) $a \vee d = 1$. By (6), we can infer that there is $P \in X(L)$ and (8) $b \vee d \not\in P$. Since there is $Q \in min \, X(L)$ such that $Q \subseteq  P$, we have by (8) that $b \vee d \not\in Q$. Therefore, $d \not\in Q$ and so, by (7) we get that $a \in Q$, which completes the proof.
\end{proof}

\

As a direct consequence of Theorem \ref{T23} and  Proposition \ref{P28} we conclude

\begin{theo}\label{T24}
Let $(L, \nabla, \triangle)$ be an m--lattice. Then the following conditions are equivalent:

\begin{itemize} 
\item[{\rm (i)}] $(L,\nabla,\triangle)$ is simple,
\item[{\rm (ii)}] $(\nabla,\triangle)$ is  simple and for each $a,\,b \in L$ if $a \not\leqslant b$, there is $c \in L$ such that $a \wedge c \not= 0$ and $b \wedge c = 0$ or there is $d \in L$ such that $b \vee d \not= 1$ and $a \vee d = 1$.
\end{itemize}
\end{theo}

Our next task is to characterize the subdirectly irreducible but not simple m--lattices. To this end, we define a new topology on $X$ whose relationship with Priestley topology is shown in Lemma \ref{L29}. 
 
\begin{lem}{\rm(\cite[Lemma 3.17]{FPZ1})}\label{L29} 
Let $(X,E)$ be an mq--space and $\tau_{\cal S}=\{X \setminus F:\, F \in {\cal C}_{idS}(X)\}$. Then, $\tau_{\cal S}$ defines a topology on $X$ whose closed sets are exactly the members of ${\cal C}_{idS}(X)$. Besides, the Priestley topology is finer than $\tau_{\cal S}$.
\end{lem}

Let $X$ be an mq--space and $Y\subseteq X$. We shall denote by $\overline Y^{\cal S}$ the closure of $Y$ when $X$ is endowed with the topology $\tau_{\cal S}$.

\begin{cor}{\rm(\cite[Corollary 3.19]{FPZ1})}\label{C29}
Let $(X,E)$ be an mq--space. Then $\overline Y = \overline Y^{\cal S}$ for all $id-$sa\-tu\-ra\-ted subset $Y$ of $X$.
\end{cor}

\begin{lem}\label{L210} 
Let $(X,E)$ be an mq--space and let $z,y \in X$ be such that $z \not\in E([y))$. Then there is $U \in D(X)$ so that $y \in \nabla_E U$ and $z \not\in \nabla_E U$. 
\end{lem}

\begin{proof} 
Since $E(z)\cap E([y))= \emptyset$ and $E([y))$ is an increasing set, we have that $t \not\leqslant w$ for all $w \in E(z)$ and  $t \in E([y))$. Then, for each $t \in E([y))$, there is $U_t \in D(X)$ such that $t \in U_t$ and $ w\not\in U_t$. Using compactness of $E([y))$ we obtain that $E([y)) \subseteq \hspace{-2mm}\bigcup\limits_{\begin{array}{c}\scriptstyle i =1 \end{array}}^{n} U_{t_i}$ further, we have $w \not\in \bigcup\limits_{ \begin{array}{c}\scriptstyle i =1 \end{array} }^{n} U_{t_i}$. Let $U_w =\hspace{-2mm}\bigcup\limits_{\begin{array}{c}\scriptstyle i =1 \end{array}}^{n} U_{t_i} $. Hence, for each $w \in E(z)$ there is $U_w \in D(X)$ such that $E([y)) \subseteq U_w$ and $w \not\in U_w$. Therefore, $E(z)\subseteq \hspace{-2mm}\bigcup\limits_{\begin{array}{c}\scriptstyle w \in E(z)\end{array}}\hspace{-2mm}(X \setminus U_{w})$ and since $E(z)$ is compact we conclude that $E(z)\subseteq \bigcup\limits_{\begin{array}{c}\scriptstyle i =1 \end{array}}^{m}(X \setminus U_{w_i})$. This last assertion implies that there is $U = \bigcap\limits_{\begin{array}{c}\scriptstyle i =1 \end{array} }^{m} U_{w_i}\in D(X)$ such that $E([y)) \subseteq U$ and $U \cap E(z)= \emptyset$, from which we conclude the proof. 
\end{proof}

\

The following results will be fundamental to prove Theorem \ref{P29}.

\begin{lem}\label{L211}
Let $(X,E)$ be a q--space. Then $E(x)$ is convex for all $x\in X$.
\end{lem}

\begin{proof} 
Let $y,z \in X$ be such that (1) $z \in E(x)$ and (2) $x \leqslant y \leqslant z$. Suppose   $(x, y)\notin E$. Then, by \cite[Lemma 2.5]{RC} there is $U \in D(X)$ such that  (3) $y \in \nabla_EU $ and (4) $x \not\in \nabla_EU $ or  (5) $x \in \nabla_EU$ and (6) $y \not\in \nabla_EU $. Assume that (3) holds. Hence, taking into account that $\nabla_EU$ is increasing and (2) we have that $z \in \nabla_EU$ and so by (1), we conclude that $x \in  \nabla_EU$ which contradicts (4). The proof of (5) is similar. Thus, we infer that $y\in E(x)$.
\end{proof}

\begin{lem}\label{L212}
Let $(X,E)$ be an  mq--space. Then there is $m \in min \,E([x))$ such that $E(x)=E(m)$. 
\end{lem}

\begin{proof}
Since $x \in E([x))$, we have by Corollary \ref{C21} that there is (1) $m \in min \,E([x))$ such that (2) $m \leqslant x$. Hence, by (1) it follows that there is $y \in X$ which verifies  $x \leqslant y$ and $(m,y)\in E$. These last assertions, (2) and  Lemma \ref{L211} imply that $(x, m)\in E$ and so, $E(x)=E(m)$.
\end{proof}

\begin{lem}\label{L213}  
Let $(X,E)$ be an mq--space and $x\in X$. If $F$ is an  $id-$saturated subset of $X$ such that $max\, E(x) \subseteq F$ then,  $max\, E(x)\cup \, min\, E([x))\subseteq F$.
\end{lem}

\begin{proof}
Let $y \in max\, E(x)$. Then, $E(x)= E(y)$ and by (mq1) we infer that  $E([x))= E([y))$.
These assertions and the fact that $F$ is an $id-$saturated subset of $X$ allow us to conclude the proof.
\end{proof}

\begin{lem}\label{L214}  
Let $(X,E)$ be an mq--space. Then for all $x \in X$, $\overline{max\, E(x)\cup}\\ \overline{\, min\, E([x))}^{\, \cal S}=  \overline{max\, E(x)}^{\, \cal S}=  \overline{min \,E([x))}^{\, \cal S}$.
\end{lem}

\begin{proof}
Since $\overline{max\, E(x)}^{\, \, \cal S}$ verifies the hypothesis of Lemma \ref{L213} we have that $max\, E(x)\cup \, min\, E([x))\subseteq  \overline{max\, E(x)}^{\, \, \cal S}$ and so, we conclude that \linebreak $\overline{max\, E(x) \cup \, min\, E([x))}^{\, \, \cal S}=  \overline{max\, E(x)}^{\, \, \cal S}$.
On the other hand, by Lemma \ref{L212}, there is (1) $m \in  min \,E([x))$ such that (2) $max\, E(x)= max\, E(m)$. From (1) and taking into account that  $\overline{min \,E([x))}^{\, \, \cal S}$  is id--saturated we have that $max\, E(m)\subseteq \overline{min \,E([x))}^{\, \, \cal S}$. Then by (2), we infer that $\overline{max\, E(x) \cup \, min \,E([x))}^{\, \, \cal S}\\ = \overline{min \,E([x))}^{\, \, \cal S}$.
\end{proof}

\begin{prop}\label{P29}
Let $(X,E)$ be an mq--space and let $Y$ be a closed subset of $X$. Then the following   conditions are equivalent:

\begin{itemize}
\item[{\rm(i)}] $Y$ is $id-$saturated,
\item[{\rm(ii)}] for all $y \in Y$, $\overline{max\, E(y)}^{\, \, \cal S}\subseteq Y$,
\item[{\rm(iii)}] for all $y \in Y$, $ \overline{min \,E([x)}^{\, \, \cal S}\subseteq Y$.
\end{itemize}
\end{prop}

\begin{proof}
(i)\, $\Rightarrow$\, (ii): From the hypothesis we have that for all $y \in Y$, $max\, E(y) \subseteq Y$. Hence, by Corollary \ref{C29} we infer that $\overline{max\, E(y)}^{\, \, \cal S}\subseteq Y$.

\vspace{2mm}

(ii)\, $\Leftrightarrow$\, (iii): It is an immediate consequence of Lemma \ref{L214}.

\vspace{2mm}

(iii)\, $\Rightarrow$\, (i): From (iii) and Lemma \ref{L214} we have that for all  $y \in Y$,\, $\overline{max\, E(y)\cup \, min\, E([y))}^{\, \cal S}\subseteq Y$ and then, we conclude the proof.
\end{proof}

\begin{theo}\label{T25}
Let $(X,E)$ be an mq--space. Then the following conditions are equivalent:
\begin{itemize}
\item[{\rm(i)}] $(D(X), \triangle_{E}, \nabla_{E})$ is a subdirectly irreducible monadic distributive lattice but not simple, 
\item[{\rm(ii)}] one and only one of these conditions hold:

\hspace{-2mm}{\rm(a)} $\{x \in X: \overline{max\, E(x)}^{\cal S}= X \}$ is a proper non--empty open subset of $X$,\\[-3mm]

\hspace{-2mm}{\rm(b)} there is $x \in X$ such that $x \not\in \overline{max\, E(x)}^{\, \cal S}$ and $X = \overline{max\,E(x)}^{\cal S}\cup \{x\}$.
\end{itemize}
\end{theo}

\begin{proof}
(i)\, $\Rightarrow$\, (ii): Let us first show that conditions (a) and (b) are incompatible. Suppose it is not the case. Then, there are elements $x, y \in X$ such that $\overline{max\, E(x)}^{\cal S}= X$,\, $y \not\in \overline{max\, E(y)}^{\cal S}$ and $X = \overline{max\, E(y)}^{\cal S} \cup \{y\}$. Then, it follows that $x \not= y$. Therefore, $x \in \overline{max\, E(y)}^{\, \cal S}$. But since $\overline{max\, E(y)}^{\, \cal S}$ is ${id-}$sa\-tu\-ra\-ted, by Proposition  \ref{P29} we infer that $\overline{max\, E(x)}^{\, \cal S} \subseteq \overline{max\, E(y)}^{\, \cal S} $ and so, $X = \overline{max\, E(y)}^{\, \cal S}$ which contradicts the fact that $y\in X \setminus \overline{max\, E(y)}^{\, \, \cal S}$.

It follows from the hypothesis and Theorem \ref{T22} that ${\cal C}_{idS}(X)\setminus \{X\}$ has a greatest element $Y$. Let $F= \{x \in X: \overline{max\, E(x)}^{\, \cal S}\not=X \}$. Since $Y$ is closed and $id-$saturated  it follows  from Proposition  \ref{P29} that for each $y \in Y$,\, $\overline{max\, E(y)}^{\, \cal S} \subseteq Y$. From this statement and taking into account that $Y\not= X$ we conclude that $\overline{max\, E(y)}^{\, \cal S} \not= X$. Therefore, $Y \subseteq F$.

Suppose first that $Y = F$. Therefore, $X \setminus F$ is a proper non--empty open subset of 
$X$ and according to the definition of $F$ we obtain (a).

Next, suppose that there is $x\in F\setminus Y$. Then, from Lemma \ref{L29} it follows that $\overline{max\, E(x)}^{\, \cal S}$ is a proper, closed and $id-$saturated  subset of $X$. Hence,\, $ \overline{max\, E(x)}^{\, \cal S}\subseteq Y$. On the other hand, $\{x\} \cup  \overline{max\, E(x)}^{\, \cal S}$ is a closed and ${id-}$saturated subset of $X$. Consequently, since $x \notin Y$ we conclude that $X = \{x\} \cup  \overline{max\, E(x)}^{\, \cal S}$. Therefore, we have proved (b).

\vspace{2mm}

(ii)\, $\Rightarrow$\, (i): Assume first that (a) holds, and $F$ be defined as above. From the hypothesis it follows that there is  $x \in F$. Besides, by Lemma \ref{L21}  we have that $max\, E(x) \not= \emptyset$. Therefore, by Lemma \ref{L29}, $\overline{max\, E(x)}^{\, \cal S}$ is a proper non--empty closed and $id-$saturated subset of $X$. This assertion implies by Theorem \ref{T22} that $D(X)$ is non--simple. On the other hand, from the hypothesis it follows that $F$ is a proper non--empty closed subset of $X$. Furthermore, $F$ is $id-$saturated. Indeed, let $x \in F$ and  $y \in \overline {max\, E(x)}^{\, \cal S}$, since  $\overline {max\, E(x)}^{\, \cal S}$ is closed and $id-$saturated we have by Proposition \ref{P29} that $\overline {max\, E(y)}^{\, \cal S} \subseteq \overline {max\, E(x)}^{\, \cal S}$.  Since $x\in F$ we conclude that $\overline {max\, E(y)}^{\, \cal S}\not= X $, which implies that $y \in F$. Moreover, take $H \in {\cal C}_{idS}(X)\setminus \{X\}$ and $h\in H$. Then,
by Proposition \ref{P29} we have that $\overline{max\, E(h)}^{\, \cal S} \subseteq H$ and since $H\not= X$ we infer that $h\in F$. Therefore, $H \subseteq F$. This means by Theorem \ref{T22}  that $D(X)$ is subdirectly irreducible.

Assume now that (b) holds. From the hypothesis it follows that $Y= \overline{max\, E(x)}^{\cal S}\in {\cal C}_{idS}(X) \setminus \{\emptyset, X\}$ and so, $D(X)$ is not a simple 
m--lattice. Let $H \in {\cal C}_{idS}(X) \setminus \{\emptyset, X\}$ and suppose that $H \not\subseteq Y$. Hence, $x \in H$ and by Proposition \ref{P29} we have that $Y \subseteq H$. Therefore, $H=X$ which is a contradiction. Then, $Y$ is the greatest element of ${\cal C}_{idS}(X)\setminus \{X\}$.
\end{proof}

\

Our next task is to give another description of subdirectly irreducible monadic distributive lattices but not simple which, in our opinion, is simplier than the one obtained above. In order to do this Propositions \ref{P210} and \ref{P211} will be fundamental. 

\begin{prop}\label{P210}
Let $(X,E)$ be an mq-space and let $x \in X$. Then the following conditions are equivalent:

\begin{itemize}
\item[{\rm(i)}] $\overline{max\, E(x)}^{\, \cal S}=X$, 
\item[{\rm(ii)}]$\overline{min \,E([x))}^{\, \cal S}=X$,
\item[{\rm(iii)}]$E([x))=X$ and $\overline{ min \,X}^{\, \cal S}= X$.
\end{itemize}
\end{prop}

\begin{proof}
(i)\, $\Leftrightarrow$\, (ii): It is a direct consequence of Lemma \ref{L214}.

\vspace{2mm}

(ii)\, $\Leftrightarrow$\, (iii): From the hypothesis and Lemma \ref{L28} we infer (iii). The reverse implication is obvious. 
\end{proof}

\begin{prop}\label{P211} 
Let $(X,E)$ be an mq--space and let $V=\{x \in X: E([x))=X \}$.  Then it holds:

\begin{itemize}
\item[{\rm(i)}] $V$ is a closed and decreasing subset of $X$ and $E(V)= V$,
\item[{\rm(ii)}] $U \in D(X)\setminus \{X\}$ and $\nabla_E(U)\not=X$ imply $\nabla_E(U) \subseteq X \setminus V$.
\end{itemize}
\end{prop}

\begin{proof}
(i)\,  Let $x \in V$ and $y \in X$ be such that $y \leqslant x$. Then, $[x) \subseteq [y)$ which implies by taking into account the definition of $V$ that $E([y))=X$ and hence, $y \in V$. Consequently $V$ is decreasing. Furthermore, let (1) $y \in X \setminus V$. Then, there is $z \in X \setminus E([y))$ and so, by Lemma \ref{L210} there is $U \in D(X)$ such that (2) $y \in \nabla_E(U)$ and  $ z \not\in \nabla_E(U)$. Therefore, $\nabla_E(U)\not=X$. From this last assertion and the fact that $E([x)) \subseteq U$ for all $x \in  \nabla_E(U)$, we conclude that $E([x)) \not= X$ hence, $\nabla_E(U) \subseteq X \setminus V$. Thus, from (1) and (2) we have that $V$ is a closed subset of $X$. Besides, $E(V)=V$. Indeed, let $y \in E(V)$. Then, there is $x \in V$ such that $E(y)= E(x)$ and by (mq1), $E([x))= E([y))$. Therefore, $E([y))=X$ and so, $y \in V$. 

\vspace{2mm}

(ii)\,  It follows from the hypothesis that $E([x)) \not=X$ for all $x \in \nabla_E(U)$ which implies that $\nabla_E(U)\subseteq X \setminus V$.
\end{proof}

\

We have now achive our desired goal.

\

\begin{theo}\label{T26}
Let $(X,E)$ be an mq--space. Then the following conditions are equivalent:

\begin{itemize}
\item[{\rm(i)}] $(D(X), \triangle_{E}, \nabla_{E})$ is a subdirectly irreducible monadic distributive lattice but not simple, 
\item[{\rm(ii)}] one and only one of these conditions hold:

\begin{itemize}
\item[{\rm(a$^{\prime}$)}] $\overline{min\, X}^{\, \cal S}= X$ and $\nabla_E(D(X))\setminus \{X\}$ has last element  which is different from $\emptyset$,
\item[{\rm(b$^{\prime}$)}] there is $x \in X$ such that $x \not\in \overline{min\, X}^{\, \cal S}$, $X =\overline{min\, X}^{\, \cal S} \cup \{x\}$ and $E([x))=X$.
\end{itemize}
\end{itemize}
\end{theo}

\begin{proof}
We only prove the equivalence between conditions (a) and (b) in Theorem \ref{T25} and 
(a$^{\prime}$) and (b$^{\prime}$) respectively.

\vspace{2mm}

(a)\, $\Leftrightarrow$\, (a$^{\prime}$): Let $W= \{x \in X:\, \overline{max\,E(x)}^{\, \cal S}= X\}$. Then, by Proposition \ref{P210} we have that (1) $\overline{min\, X}^{\, \cal S}= X$ and  $W= \{x \in X:\, E([x))=X\}$. This last assertion and condition (i) in Proposition \ref{P211} imply that $W$ is a closed decreasing subset of $X$ and $E(W)=W$. Therefore, from (a) we infer that $X \setminus W \in \nabla_E(D(X))\setminus \{\emptyset, X\}$. Besides, (ii) in   
Proposition \ref{P211} allows us to conclude that (2) $X \setminus W$ is the last element of 
$\nabla_E(D(X))\setminus \{X\}$ and $X \setminus W \not=\emptyset$. Hence, from (1) and (2) we obtain (a$^{\prime}$). Conversely, let (3) $U$ be the last element of $\nabla _E(D(X))\setminus \{ X\}$, $U\not=\emptyset$. Hence, $\nabla _E(U)\not=X$ and so, by (ii) in Proposition \ref{P211} we have that $U \subseteq X \setminus V$. Suppose now  that $U \subset X \setminus V$. Then, there is $x \in (X \setminus V)\cap(X \setminus U)$ and bearing in mind the definition of $V$ we conclude that $E([x)) \not=X$. This means that there is $y \not\in E([x))$ and hence, by Lemma \ref{L210} it follows that there is  $W \in D(X)$ such that  
$x \in \nabla_E(W)$ and $y \not\in \nabla_E(W)$. Hence, $\nabla_E(W) \in \nabla_E(D(X))\setminus \{X, \emptyset\}$ and  $\nabla_E(W) \not\subseteq U$ which contradicts (3). Therefore, (4) $U = X \setminus V$. On the other hand, since $\overline{min\, X}^{\, \cal S}= X$ we infer from Proposition \ref{P210} that $V= \{x \in X: \overline{max\,E(x)}^{\, \cal S}= X\}$ and so, by (3) and (4) the proof of (a) is complete.  

\vspace{2mm}

(b)\, $\Leftrightarrow$\, (b$^{\prime}$):  
From  (i) in Lemma \ref{L28} and Proposition \ref{P29}, it follows that $\overline{max\, E(x))}^{\, \cal S}\subseteq E([x))$. Since $X=\overline{max\, E(x))}^{\, \cal S} \cup \{x\}$  we conclude that (1) $E([x))=X$. On the other hand, from the hypothesis, Lemma \ref{L214} and (1) we infer that $x \not\in \overline{min\, X}^{\, \cal S}$ and $X=\overline{min\, X}^{\, \cal S}\cup \{x\}$. Hence, we have shown (b$^{\prime}$). Conversely, since $X = E([x))$ from   Lemma \ref{L214} we obtain that $\overline{min\, X}^{\, \cal S}= \overline{max\, E(x)}^{\, \cal S}$ and so, we conclude (b). 
\end{proof}

\begin{cor}\label{C210} 
Let $(X,E)$ be an mq--space and let $(\triangle_{E}, \nabla_{E})$ be simple. Then the following conditions are equivalent:

\begin{itemize}
\item[{\rm(i)}] $(D(X),\triangle_{E}, \nabla_{E})$ is a subdirectly irreducible monadic distributive lattice but not simple, 

\item[{\rm(ii)}] there is $x \in X$ such that $x \not\in \overline{min\, X \cup \, max\, X}$,  
 $X =  \overline{ min\, X \cup \, max\, X } \cup \{x\} $.
\end{itemize}
\end{cor}

\begin{proof}
(i)\, $\Rightarrow$\, (ii):  Since $(\triangle_{E}, \nabla_{E})$ is simple, then $\nabla_E(D(X))\setminus \{X\}=\{\emptyset\}$. Therefore, condition (a$^{\prime}$) in Theorem \ref{T26} is not verified from which it follows that condition (b$^{\prime}$) holds. Hence, taking into account Corollary \ref{C28} we have that $\overline{min\, X}^{\, \cal S} = \overline{min\, X \cup \, max\, X}$ and so, we conclude the proof.

\vspace{2mm}

(ii)\, $\Rightarrow$\, (i): Proposition \ref{P26} and the fact that $(\triangle_{E}, \nabla_{E})$ is simple imply that $E= X \times X$ and therefore, $E([x))=X$ for all $x \in X$. Besides, 
from (ii) and  Corollary \ref{C28} we have that $x \not\in \overline{min\, X}^{\, \cal S}$ and   $X = \overline{min\, X}^{\, \cal S} \cup \{x\}$. Then, by (b$^{\prime}$) in Theorem \ref{T26} the proof is complete.
\end{proof}

\section{\large Monadic augmented Kripke frames}\label{S3}\indent

Our next task is to show the relationship between the categories $\mbox{\boldmath ${\cal P K F}$}$ and $\mbox{\boldmath $m{\cal Q}$}$. To this end, we determine in the first place a new  topological duality for monadic distributive lattices by considering the category whose objects are augmented Kripke frames which verify certain additional conditions. More precisely,

\begin{defi}\label{D31}
A monadic augmented Kripke frame {\rm(}or mk--frame{\rm)} is a quadruple $(X, \Omega, \leqslant, E)$ where $(X,\leqslant)$ is a non--empty partially ordered set, $E$ is an equivalence relation on $X$ and the following conditions are verified:

\begin{itemize}
\item[\rm{(mk1)}]$(X, \leqslant, E)$ is an augmented Kripke frame,
\item[\rm{(mk2)}]$(X, \Omega, \leqslant)$ is a Priestley space,
\item[\rm{(mk3)}]$E$ is a closed relation,
\item[\rm{(mk4)}] for each $U \in D(X)$, $E(U)$ is an open subset of $X$,
\item[\rm{(mk5)}] for each $U \in D(X)$, $(E(X \setminus U)]$ is an open subset of $X$.
\end{itemize}
\end{defi} 

In what follows, we will denote monadic augmented Kripke frames by $(X,\leqslant, E)$.

\begin{defi}\label{D32} 
Let $(X_1,\leq_1, E_1)$ and $(X_2,\leq_2, E_2)$ be mk--frames. An mk--function $f: X_1 \fun X_2$ is an order--preserving continuous function  which verifies the following conditions:

\begin{itemize}
\item[\rm{(mkf1)}] $(x,y)\in E_1$ implies $(f(x), f(y)) \in E_2$,
\item[\rm{(mkf2)}] $E_2(f(x))\subseteq  (f(E_1(x))]$ for all $x \in X_1$,
\item[\rm{(mkf3)}] $E_2([f(x)) )\subseteq  [f(E_1(x)))$ for all $x \in X_1$.
\end{itemize}
\end{defi}

The category of mk-frames and mk--functions will denote by $\mbox{\boldmath $m{\cal K F}$}$.
		
\begin{prop}\label{P31} 
Let $X$ be a non--empty set. Then the following conditions are equivalent:

\begin{itemize}
\item[{\rm(i)}] $(X, \Omega, \leqslant, E)$ is an mq--space,
\item[{\rm(ii)}] $(X, \Omega, \leqslant, E)$ is an mk--frame.
\end{itemize}
\end{prop}

\begin{proof}
(i)\, $\Rightarrow $\, (ii): From Remark \ref{R21} it follows that $(X,\Omega, \leqslant, E)$ is an augmented Kripke frame. Then, Lemma \ref{L21} allows us to complete the proof.

\vspace{2mm}
(ii)\, $\Rightarrow $\, (i): Since $X$ is a  Hausdorff space, for each $x \in X$ we have that $\{x\}$ is closed, then by (mk3) it follows that $E(x)$ is a closed  subset of $X$. On the other hand, from (mk3) and (mk4) we get that $E(U)$ is clopen for all $U \in D(X)$. Besides,  $E(U)$ is an increasing subset of $X$. Indeed, let $y \in E(U)$ and $y \leqslant z$. Then, there is (1) $x \in U$ such that $y \in E(x)$. Hence, $z \in [E(x))$ and from (mk1) we conclude that $z \in E([x))$. Furthermore, from (1) we infer that $E([x))\subseteq E(U)$ and so, $z \in E(U)$. Therefore, $(X, \Omega, \leqslant, E)$ is a q--space. Moreover, by (mk3) we have that for each $V \in D(X)$, $E(X \setminus U)$ is a closed subset of $X$ and since $X$ is a Priestley  space we conclude that $(E(X \setminus U)]$ is closed. From this last assertion and (mk5) we obtain (mq2) and so, the proof is complete. 
\end{proof}

\

Next, we are going to show that the notions of mq--function and mk--function are also equivalent. To this end, first we will indicate a cha\-rac\-te\-ri\-za\-tion of q--functions  proved in \cite{FPZ}, from which we obtain a new description of mq--functions.

\begin{prop} \label{P32}{\rm(\cite[Proposition 2.1]{FPZ})} 
Let $(X_1, E_1)$ and $(X_2, E_2)$ be q--spaces and  $f$ an order--preserving continuous function from $X_1$ into $X_2$. Then the following conditions are equivalent:\vspace{-3mm}

\begin{itemize}
\item[{\rm(i)}] $f$ is a q--function, 
\item[{\rm(ii)}] $f$ satisfies 

\begin{itemize}
\item[{\rm(f1)}] $(x,y)\in E_1$ implies $(f(x),f(y))\in E_2$,
\item[{\rm(f2)}] $E_2(f(x))\subseteq (f(E_1(x))]$  for all $ x \in X_1$.
\end{itemize}
\end{itemize}
\end{prop}

\begin{prop}\label{P33}
Let $(X_1, E_1)$, $(X_2, E_2)$ be mq--spaces and $f:X_1\fun X_2$\, a function. Then the following conditions are equivalent:

\begin{itemize}
\item[\rm (i)] $f$ is an mq-function,
\item[\rm (ii)] $f$ is an order--preserving continuous function satisfying {\rm (f1)}, {\rm (f2)} and  {\rm (mqf3)}. 
\end{itemize}
\end{prop}

\begin{proof} 
It is a direct consequence of Proposition \ref{P32} and  Corollary \ref{C22}.
\end{proof}

\begin{cor}\label{C31}
Let $(X_1,E_1)$, $(X_2,E_2)$ be mq--spaces and $f:X_1\fun X_2$\, a function. Then the following conditions are equivalent:

\begin{itemize}
\item[\rm (i)] $f$ is an mq--function,
\item[\rm (ii)] $f$ is an mk--function. 
\end{itemize}
\end{cor}

\begin{proof}
It follows as  an immediate consequence of Proposition \ref{P33}.
\end{proof}

\

Theorem \ref{T21}, Proposition \ref{P31} and  Corollary \ref{C31} allow us to conclude Theorem \ref{T31}.

\begin{theo}\label{T31}
The categories $\mbox{\boldmath ${\cal M}$}$ and $\mbox{\boldmath $m{\cal K F}$}$ are  dually e\-qui\-va\-len\-t.
\end{theo}

From now until the end of this section, our attention  is focused on determining the relationship between $\mbox{\boldmath ${\cal P K F}$}$ and $\mbox{\boldmath $m{\cal K F}$}$.

\begin{prop}\label{P34}
Let $(X, \Omega,R)$ be a perfect Kripke frame. If $(X,R)$ is a partially ordered set then, $(X, \Omega,R)$ is a Priestley space.
\end{prop}

\begin{proof} 
From the hypothesis $(X,\Omega)$ is a compact and Hausdorff topological space. Then, we only need to show that $(X,\Omega,R)$ is totally order--disconnected topological space. Let $x, y \in X$ be such that (1) $y \not\in R(x)$. Since $R(x)$ is a closed set and the set ${\cal CP}(X)$ of all clopen subsets of $X$ is a base for $\Omega$, there is $A \in {\cal CP}(X)$ such that $y \in A$ and $A \cap R(x) = \emptyset$. Hence, (2) $x \not\in R^{-1}(A)$ and from the fact that $R$ is a reflexive relation we have that (3) $y \in R^{-1}(A)$. On the other hand,  $R$ is a perfect relation which implies that $R^{-1}(A)\in {\cal CP}(X)$ and so, (4) $X \setminus  R^{-1}(A)\in {\cal CP}(X)$.  Besides, $X \setminus  R^{-1}(A)$ is an increasing subset of $X$. Indeed, let $z \in X \setminus  R^{-1}(A)$ and (5) $w \in R(z)$. Then, we have that (6) $R(z) \cap A = \emptyset$ and taking into account that $R$ is a transitive relation from (5) it follows that (7) $R(w)\subseteq R(z)$. Hence, from (6) and (7) we infer that $R(w) \cap A = \emptyset$ and so, $w \in X \setminus  R^{-1}(A)$. Therefore,
$X \setminus  R^{-1}(A) \in D(X)$ and from (2) and (3) we conclude the proof.
\end{proof}

\begin{cor}\label{C32} 
Every perfect augmented Kripke frame is a monadic augmented Kripke frame.
\end{cor}

\begin{proof}
It is a direct consequence of Proposition \ref{P34} and Lemma 7 in \cite[Section 4]{GBz}.
\end{proof}

\

In general, the converse of Corollary \ref{C32} is not true as the following example shows:

\vspace{2mm}

\begin{ex}\label{E31} 
Let $\RR$ be the set of real numbers endowed by the Euclidean topology and $\cal F$ the set of all closed subsets of $\RR$. It is well known that $\langle {\cal F},\cap , \cup,\emptyset, \RR \rangle$ is a bounded distributive lattice. Besides, $({\cal F},\nabla,\triangle)$ is a monadic distributive lattice where the operators $\nabla,\triangle$ are defining by the prescriptions $\nabla \emptyset=\emptyset$ and $\nabla F=\RR$ for each $F\not=\emptyset$; $\triangle \RR=\RR$ and $\triangle F=\emptyset$ for each $F\not=\RR$. 
Then, the monadic augmented Kripke frame $(X({\cal F}),\subseteq, E_\nabla)$ associated
with ${\cal F}$ is not perfect since for all $U \in D(X({\cal F}))$ we have that $(U]$ is not  an open subset of $X({\cal F})$. Indeed, if it were, it follows that $X({\cal F})\setminus (U] \in D((X({\cal F}))$ for all  $U \in D(X({\cal F}))$ and  therefore, $F\to \emptyset$ would be defined, which is a contradiction.
\end{ex}

\

Our next task will be to show that the morphisms between perfect Kripke frames are also morphisms between monadic augmented Kripke frames. First, we will determine properties of    
mk-frames which will be useful to this aim. 

\begin{lem}\label{L31} 
If $(X,\leqslant,E)$ is an mk--frame, then

\begin{itemize}
\item[\rm (i)] for each $x\in X$, $max\, E_{E \circ \leqslant}(x)\not= \emptyset$,
\item[\rm (ii)] $E = E_{E \circ \leqslant}$. 
\end{itemize}
\end{lem}

\begin{proof}
(i)\, From (mk2) and (mk3) it follows that for each $x\in X$, $E([x))$ and $(E(x)]$ are closed subsets of $X$. Then, taking into account that $E_{E \circ \leqslant}(x) = E([x))\cap (E(x)]$ we conclude the proof.

\vspace{2mm}

(ii)\, It follows from (mk1), (i) and  Lemma 3 in \cite[Section 2]{GBz}.
\end{proof}

\

\begin{rem}\label{R31}
Note that strongly isotone maps are also isotone.
\end{rem}

\begin{lem}\label{L32}
Let $(X_1,\leqslant_1,E_1)$, $(X_2, \leqslant_2, E_2)$ be mk--frames and let $f: X_1\fun X_2$ be strongly isotone with respect to $E_1\circ \leqslant_1$. Then,

\begin{itemize}
\item[{\rm(i)}] $f$ is isotone with respect to $E_1$,
\item[{\rm(ii)}] $f$ is almost strongly isotone with respect to $E_1\circ \leqslant_1$.
\end{itemize}
\end{lem}

\begin{proof} 
(i)\, Let $(x,y)\in E_1$. Then, by Lemma \ref{L31} we have that $(x,y)\in E_1\circ 
\leqslant_1$ and  $(y,x)\in E_1\circ \leqslant_1$. From the hypothesis and Remark \ref{R31}  we conclude that $(f(x),f(y))\in E_2 \circ \leqslant_2 $ and $(f(y),f(x))\in E_2\circ \leqslant_2$. Therefore, $(f(x),f(y))\in E_{E_2\circ \leqslant_2}$ and so, by 
Lemma \ref{L31} we have that $(f(x),f(y))\in E_2$.

\vspace{2mm}

(ii)\, It is straightforward. 
\end{proof}

\

From Remark \ref{R31} and Lemma \ref{L32} we have the following statement.

\begin{prop}\label{P35}
Every morphism in $\mbox{\boldmath ${\cal P A K F}$}$ is a morphim in $\mbox{\boldmath $m{\cal  K F}$}$.
\end{prop}

\vspace{2mm}

Corollary \ref{C32}, Example \ref{E31} and Proposition \ref{P35} allows us to conclude 

\begin {theo}\label{T32} 
The category $\mbox{\boldmath ${\cal P A K F}$}$ is a proper subcategory of $\mbox{\boldmath $m{\cal  K F}$}$.
\end{theo}

\

\noindent  A. V. Figallo, Departamento de Matem\'atica, Universidad Nacional del Sur, 
8000 - Bah\'{\i}a Blanca, Argentina.\\
\noindent Instituto de Ciencias B\'asicas. Universidad Nacional de  San Juan. 5400 - San Juan, Argentina.{\it e-mail:avfigallo@gmail.com}\\[2mm]

\noindent I. Pascual, Instituto de Ciencias B\'asicas. Universidad Nacional de San Juan, 5400 - San Juan, Argentina.{\it e-mail: inespascual756@hotmail.com}\\[2mm]

\noindent Alicia N. Ziliani, Departamento de Matem\'atica, Universidad Nacional del Sur, 8000 - Bah\'{\i}a Blanca, Argentina.\\ 
\noindent Instituto de Ciencias B\'asicas. Universidad Nacional de  San Juan. 5400 - San Juan, Argentina. {\it e-mail: aziliani@criba.edu.ar}

\begin{thebibliography}{99} 

\bibitem{RB.PD} R. Balbes and P. Dwinger, {\em Distributive lattices}, University of Mis\-sou\-ri Press, Columbia, 1974.


\bibitem{GBz0} G. Bezhanishvili, {\em Varieties of monadic Heyting algebras. Part I}, Studia Logica 61, 3(1998), 367--402.


\bibitem{GBz} G. Bezhanishvili, {\em Varieties of monadic Heyting algebras. Part II: Duality Theory}, Studia Logica 62(1998), 1--28. 

\bibitem{GBz1} G. Bezhanishvili, {\em Varieties of monadic Heyting algebras. Part III}, Studia Logica 64, 2(2000), 215--256. 

\bibitem{GB}  G. Birkhoff, {\em  Lattice theory}, Amer. Math. Soc., Col Pub., 25 3rd ed., Providence, 1967. 

\bibitem{RC} R. Cignoli, {\em Quantifiers on distributive lattices}, Discrete Math., 96(1991), 183--197.

\bibitem{FZ} A. V. Figallo and A. Ziliani, {\em Notes on monadic distributive lattices}, Preprints del Instituto de Ciencias B\'asicas, U. N. de San Juan, Argentina, 2, 1(1997), 19--35.

\bibitem{FPZ} A. V. Figallo, I. Pascual and A. Ziliani, {\em Notes on monadic $n$-valued Lukasiewicz  algebras}, Math. Bohemica 129, 3(2004), 255--271. 

\bibitem{FPZ1} A. Figallo, I. Pascual and A. Ziliani, {\em Monadic distributive lattices},
Logic Jnl IGPL, 15(2007), 535--551. 

\bibitem{MacL} S. Mac Lane, {\em Categories for the Working Mathematician}, Springer--Verlag, Berlin, 1971.

\bibitem{AM.OV} A. Monteiro and O. Varsavsky, {\em Algebras de Heyting mon\'adicas}, Actas de las X Jornadas de la Uni\'on Matem\'atica Argentina, Bah\'{\i}a Blanca, (1957), 52--62. 
(A French translation is published as Notas de L\'ogica Matem\'atica 1, Instituto de Matem\'atica, Universidad Nacional del Sur, Bah\'{\i}a Blanca (1974), 1--16.)  

\bibitem{HP1} H. Priestley, {\em Representation of distributive lattices by means of ordered Stone spaces}, Bull. London Math. Soc., 2(1970), 186--190.

\bibitem{HP2} H. Priestley, {\em Ordered topological spaces and the representation of dis\-tri\-bu\-ti\-ve lattices}, Proc. London Math.  Soc., 4, 3(1972), 507--530.

\bibitem{HP3} H. Priestley, {\em Ordered sets and duality for distributive lattices}, Ann. Discrete Math., 23(1984), 39--60.
\end{thebibliography}
\end{document}